\documentclass[12pt,a4paper]{article}
\usepackage{amsmath,amssymb,amsfonts}
\def\Bbb{\mathbb}

\title{\bf On the distribution of Dedekind sums}

\author{Kurt Girstmair}
\date{}

\makeatletter
\let\@@maketitle=\maketitle
\def\maketitle{\def\thispagestyle##1{\relax}\@@maketitle}
\makeatother
\def\BE{\begin{equation}}
\def\EE{\end{equation}}
\def\BD{\begin{displaymath}}
\def\ED{\end{displaymath}}
\def\BA{\begin{array}}
\def\EA{\end{array}}
\def\BEA{\begin{eqnarray*}}
\def\EEA{\end{eqnarray*}}
\def\BI{\bibitem}

\def\Z{\Bbb Z}
\def\Q{\Bbb Q}
\def\R{\Bbb R}

\def\FF{{\cal F}}

\def\phi{\varphi}
\def\EPS{\varepsilon}

\def\MB{\mbox}
\def\LD{\ldots}
\def\OV{\overline}

\def\sminus{\smallsetminus}
\def\DIV{\,|\,}
\def\NDIV{\, \nmid \,}

\def\BQ{``}
\def\EQ{'' }

\def\MN{\medskip\noindent}

\def\DED{Dedekind }

\def\PKT{\circle{0.06}}
\def\PT#1#2{\put(#1,#2){\PKT}}

\begin{document}
\maketitle

\begin{abstract}
\noindent
Dedekind sums have applications in quite a number of fields of mathematics. Therefore, their distribution has found considerable interest.
This article gives a survey of several aspects of the distribution of these sums. In particular, it highlights results
about the values of Dedekind sums, their density and uniform distribution. Further topics include mean values,
large and small (absolute) values, and the behaviour of Dedekind sums near quadratic irrationals.
The present paper can be considered as a supplement to the survey article [R. W. Bruggeman, On the distribution of Dedekind sums,
Contemp. Math. 166 (1994), 197--210].
\end{abstract}

\section*{1. Introduction}

Let $b$ be a positive integer and $a\in \Z$, $(a,b)=1$. The classical \DED sum $s(a,b)$ is defined by
\BD
   s(a,b)=\sum_{k=1}^b ((k/b))((ak/b))
\ED
where $((\LD))$ is the \BQ sawtooth function\EQ defined by
\BE
\label{1.1}
  ((t))=\begin{cases}
                 t-\lfloor t\rfloor-1/2, & \MB{ if } t\in\R\smallsetminus \Z; \\
                 0, & \MB{ if } t\in \Z.
               \end{cases}
\EE
(see, for instance, \cite[p. 1]{RaGr}). In many cases it is more
convenient to work with
\BD
 S(a,b)=12s(a,b)
\ED instead.
We call $S(a,b)$ a {\em normalized} \DED sum.

\DED sums have many interesting applications, for instance, in
the theory of modular forms (see \cite{Ap, Br}), in algebraic number theory (class numbers, see \cite{BeEv, Me}),
in connection with lattice point problems (see \cite{BeRo, RaGr}), topology (see \cite{At, KiMe}) and algebraic geometry (see \cite{Ur}).
Various generalizations of \DED sums have been introduced for similar purposes (see \cite{Sz, Za}). In \cite{Kn}
such a generalization is used for the assessment of random number generators.
We emphasize that this enumeration is by no means exhaustive (nor are the corresponding references).

The general interest in \DED sums justifies the study of the distribution of these rational numbers.
Several authors have contributed to this subject (see, for instance, \cite{Br, CoFr, Gi, Hi,  Va, Zh}).

The present article is meant as a survey of the distribution of \DED sums --- we do not pretend that this survey is complete.
In particular, we consider
\begin{itemize}
\item the set of values of \DED sums,
\item the density of this set and the question of uniform distribution,
\item mean values,
\item large and small \DED sums,
\item the behaviour of \DED sums near quadratic irrationals.
\end{itemize}

In addition to proved results, we also present a number of conjectures and open questions.
We should like to draw the reader's attention to the survey article \cite{Br}, most of whose results are {\em not} rendered here.

\section*{2. Some basic tools}

Given $a, b\in\Z$, $b\ge 1$, $(a,b)=1$, we call $b$ the {\em parameter} and $a$ the {\em argument} of the normalized \DED sum $S(a,b)$.
From the definition (\ref {1.1}) one sees
\BE
\label{2.1}
S(a+b,b)=S(a,b) \MB{ and } S(-a,b)=-S(a,b).
\EE
Accordingly, \DED sums $S(a,b)$ are periodic mod $b$ and odd functions in their arguments $a$. In particular, we obtain all possible values $S(a,b)$
for a fixed parameter $b$, if we restrict $a$ to the range $0\le a<b$.

Probably the most frequently used property of normalized \DED sums is the {\em reciprocity law} (see \cite[p. 3]{RaGr}):
If $a, b$ are positive integers, then
\BD
  S(a,b)+S(b,a)=\frac ab+\frac ba+\frac 1{ab} -3.
\ED

Many of the results presented here have been proved by means of the {\em three-term relation}:
Let $a,b,c,d$ be integers, $b,d>0$, $(a,b)=(c,d)=1$, and define $q=ad-bc$ and $r=aj-bk$,
where $j,k$ are such that $-cj+dk=1$. Let $q\ne 0$. Then
\BE
\label{2.3}
  S(a,b)=S(c,d)+\EPS\cdot S(r,|q|)+\frac b{dq}+\frac d{bq}+\frac q{bd}-3\EPS,
\EE
where $\EPS$ is the sign of $q$ (see, for instance, \cite[formula (9)]{Gi}).

Another important tool is the connection of normalized \DED sums with continued fractions, as expressed by the{\em  Barkan-Hickerson-Knuth formula} (see, for example, \cite{Hi}).
Let $a,b$ be integers, $0< a<b$, $(a,b)=1$. Suppose that $a/b$ equals the regular continued fraction $[0,c_1,\LD,c_n]$.
Then
\BE
\label{2.5}
  S(a,b)=\sum_{j=1}^n(-1)^{j-1}c_j+\frac{a+a^*}b+\begin{cases} -3, & \MB{ if } n \MB{ is odd;}\\
                                                               -1, & \MB{ otherwise.}
                                                 \end{cases}
\EE
Here $a^*$ is defined by $0<a^*<b$ and $aa^*\equiv 1\mod b$.

Some results mentioned here use connections of \DED sums with modular forms (see \cite{Br}) or values of $L$-series (see \cite{Zh}). Here we do not go into details.
\section*{3. The values of \DED sums}

It is known that $bS(a,b)$ is an integer (see \cite[p. 27]{RaGr}).
Accordingly, $S(a,b)$ is a rational number $k/q$, $k, q\in\Z$, $q\ge 1$, $(k,q)=1$, such that $q$ divides the parameter $b$.
So far it is not known
which numerators $k$ are possible for a given denominator $q$ of a normalized \DED sum (a problem already mentioned in \cite[p. 28]{RaGr}). However, it has been shown in \cite{Gi2}
that these numerators form {\em complete} residue classes mod $q(q^2-1)$. In other words, a value $k/q$ does not appear isolated, but all
numbers $k/q+r(q^2-1)$, $r\in\Z$, are also values of normalized \DED sums. The numerators $k$ are subject to the following congruence
conditions.

(a) If  $3\NDIV q$, then $k\equiv 0 \mod 3$.

(b) If $2\NDIV q$, then
\BD
   k\equiv\begin{cases} 2\mod 4,& \MB{ if } q\equiv 3\mod 4;\\
                        0\mod 8,& \MB{ if } q \MB{ is a square;}\\
                        0\mod 4,& \MB{ otherwise.}
          \end{cases}
\ED
Many residue classes (mod $q(q^2-1)$) of numbers $k$ that occur as numerators can be found by means of a search procedure described in \cite{Gi2}.
In this way it could be shown, for $1\le q\le 60$, that all numbers $k$, $(k,q)=1$, satisfying the conditions (a) and (b)
belong to actual values $k/q$ of normalized \DED sums. The main tool used in \cite{Gi2} is the three-term relation (\ref{2.3}).

Hence one may conjecture that these necessary conditions, together with $(k,q)=1$,
are sufficient for $k/q$ being the value of a normalized \DED sum. It would be highly desirable to prove this conjecture or, if it is false,
to extend the conditions (a) and (b) to a set of sufficient conditions.

Let us briefly look at the case $q=60$. Here conditions (a) and (b) are empty, so one expects $\phi(q)(q^2-1)=57584$ residue classes of numerators $k$.
Because of the second identity of (\ref{2.1}), it suffices to find one of the classes $\OV k$, $\OV{-k}$. Accordingly, the search procedure in question had to exhibit
$28792$ residue classes. For this purpose it needed about $2.5$ million pairs $(a,b)$.

For each value $k/q$, $(k,q)=1$, of a normalized \DED sum there are infinitely many parameters $b$ such that $S(a,b)=k/q$ for some argument $a$ (see \cite{Gi3}).

Several authors have discussed another type of values, namely, the values of the integers $bS(a,b)$ (see \cite{As, My, Na, Sai, Sal}).
For instance, it was shown that $\pm 24$, $\pm 34$ and $\pm 88$ do not have the form $bS(a,b)$ for any choice of $a$ and $b$ (see \cite{As}).
This list of exceptional values was completed in \cite{Sai}.

In our opinion, however, these results rather concern possible {\em parameters} $b$ than {\em values} of $S(a,b)$. We give an example.  Suppose that $bS(a,b)=bk/q=24$.
Then \cite[Th. 3]{Gi2} shows that $k$ must be either $24$ or $12$.
The case $k=24$ yields $b=q$, and $k=12$ yields $b=2q$. However, the number 24 is {\em not} of the form $bS(a,b)$. Accordingly, the cases $k=24, b=q$ and $k=12, b=2q$,
are impossible. For example, if $k=24$ and $q=5$, we obtain $k/q=S(3,25)$. So $b=25$ is a parameter that yields $S(a,b)=24/5$ for $a=3$. But $b=q=5$ is not a parameter of this kind, since $S(a,5)$ takes only the values
$0$ and $\pm 12/5$. In the same way $S(1,5)=12/5$. Hence $5$ is a possible parameter for $k=12$ and $q=5$, whereas $2q=10$ is not, since $S(a,10)\in\{0, \pm 36/5\}$ for the respective arguments $a$.

Whereas the value of $S(a,b)$ is subject to the restrictions (a) and (b), no restrictions occur for the {\em fractional part} of $S(a,b)$. In other words, every rational number $r$, $0\le r<1$,
is the fractional part of some normalized \DED sum, see \cite{Gi4}.

\section*{4. Density and uniform distribution}

In \cite{Hi} it has been shown that the set
\BD
\{(a/b,S(a,b)); a,b\in \Z, b\ge 1, (a,b)=1\}
\ED
is dense in the plane $\R^2$.  In particular,
\BD
\{S(a,b);a,b\in\Z, b\ge 1, 0\le a<b, (a,b)=1\}
\ED
is dense in $\R$. The main tools of \cite{Hi} are continud fractions, in particular, formula (\ref{2.5}). Let $x\in\Q$ and $\EPS>0$ be given. The paper \cite{Gi5} explicitly describes $a$ and $b$ (in terms of $x$ and $\EPS$) such that
$|S(a,b)-x|<\EPS$. A similar result is contained in the forthcoming article \cite{Ko}, which, however, leads to a smaller value of the parameter $b$
than that of \cite{Gi5}.

The density of the set of normalized \DED sums in the $p$-adic number field $\Q_p$ has been investigated in \cite{Gi6}. In the case of $p\in\{2,3\}$, normalized
\DED sums do not approximate $p$-adic units, so they are not dense in $\Q_p$. For $p\ge 5$, they are dense in $\Q_p$.

A quite different question concerns the {\em density of the set of parameters} $b$ for a given value $k/q$ of normalized \DED sums. The sequence of parameters $b$ that
was exhibited in \cite{Gi3} grows exponentially --- so it is rather thin within the set of positive integers. In a number of cases, the author of this article could give
such a sequence that grows like $Cn^4$, $n$ running through the set of positive integers (still unpublished). Hence this question waits for further investigations.

Results about the uniform distribution of \DED sums can be found in \cite{BaSh, Br2, Br, My2, Va2}.
We only mention the following result of \cite{Br}. Let $f$ be a real-valued, continuous function on $\R/\Z\times \R$ with
compact support. For $X>0$, define
\BD
  U_X(f)=\sum_{1\le b<X}\sum_{\genfrac{}{}{0pt}{1}{0\le a<b,}{(a,b)=1}} f(a/b,S(a,b)).
\ED
Then
\BE
\label{4.3}
 U_X(f)\sim\frac{X^2}{\log X}\cdot\frac{1}{2\pi^2}\int_{\R/\Z}\int_{\R}f(x,y)dy\,dx
\EE
as $X$ tends to infinity (asymptotic equality). Like other distribution results of \cite{Br}, this assertion is proved by means of Fourier coefficients of real analytic
modular forms. The reader may also look at the diagrams in \cite{Br}, which illustrate these results.

\section*{5. Mean values}

The properties (\ref{2.1}) of normalized \DED sums imply
\BD
   \sum_{\genfrac{}{}{0pt}{1}{0\le a<b,}{(a,b)=1}}S(a,b)=0.
\ED
Hence only the arithmetic mean of $|S(a,b)|$, for varying arguments $a$, is of interest. Here some insight is given by formula (\ref{2.5}) and mean value results for continued
fractions (see \cite{GiSch}).  In this way one obtains
\BE
\label{4.1}
  \frac{1}{\phi(b)}\sum_{\genfrac{}{}{0pt}{1}{0\le a<b,}{(a,b)=1}}|S(a,b)|\le \frac{6}{\pi^2}\log^2b+O(\log b)
\EE
as $b$ tends to infinity. In the said paper it is shown that
\BD
  \frac{1}{\phi(b)}\sum_{\genfrac{}{}{0pt}{1}{0\le a<b,}{(a,b)=1}}|S(a,b)|\ge \frac{3}{\pi^2}\log^2b+O(\log^2 b/\log\log b)
\ED
for $b\to\infty$. Computations suggest, however, that (\ref{4.1}) remains true if \BQ$\le$\EQ is replaced by \BQ$=$\EQ and the error term on the right hand side by a less good one.
This has to do with the fact that the summation process of \cite{GiSch} (which is based on (\ref{2.3})) can consider only \BQ large\EQ values $|S(a,b)|$.
Hence it would be desirable to show
\BD
  \frac{1}{\phi(b)}\sum_{\genfrac{}{}{0pt}{1}{0\le a<b,}{(a,b)=1}}|S(a,b)|\sim \frac{6}{\pi^2}\log^2b
\ED
for $b\to\infty$.

The quadratic mean value of $S(a,b)$ for varying arguments $a$ has been determined in \cite{CoFr}. A sharper result of \cite{Zh} says
\BD
\frac{1}{\phi(b)}\sum_{\genfrac{}{}{0pt}{1}{0\le a<b,}{(a,b)=1}}|S(a,b)|^2=5\lambda(b)\cdot b+O\left(\exp\left(\frac{4\log b}{\log\log b}\right)\cdot \frac{b}{\phi(b)}\right).
\ED
with
\BD
\lambda(b)=\frac{\prod_{p^k\,||\, b}((1+1/p)^2-1/p^{3k+1})}{\prod_{p\DIV b}(1+1/p+1/p^2)}.
\ED
As usual, $p$ runs through the prime divisors of $b$ and $p^k\,||\,b$ means that $p^k$ is the largest power of $p$ dividing $b$.
The quantity $\lambda(b)$ grows at most like $\log\log b$. Hence $\sqrt b$ is roughly the order of magnitude of the quadratic mean value of the normalized \DED sums $S(a,b)$.

Higher power mean values can be found in \cite{CoFr}.

\section*{6. Large and small values}

First observe that $|S(a,b)|\le S(1,b)<b$ for all arguments $a$ with $(a,b)=1$ (see \cite[Satz 2]{Ra}).
The mean values of the foregoing section give an idea of what \BQ large\EQ and \BQ small\EQ stand for. The quadratic mean value $\approx \sqrt b$ of the $S(a,b)$ for varying
arguments $a$ suggests that \BQ large\EQ means an order of magnitude $\gg\sqrt b$ for $|S(a,b)|$, or, in a somewhat wider sense, $\gg b^{\alpha}$ for some $\alpha$, $0<\alpha\le 1$,
if $b$ tends to infinity.
On the other hand, the arithmetic mean of the \DED sums suggests that in general $|S(a,b)|\ll \log^2 b$, so \BQ small\EQ refers to \DED sums of logarithmic size.

The first result describes a subset of the interval $I=[0,b]$ outside which all \DED sums are $\le 3\sqrt b +5$. To this end let $d$ run through all integers $1\le d\le\sqrt b$,
and $c$ through the integers $0\le c\le d$, $(c,d)=1$. Define the interval
\BD
   I_{c/d}=\{x\in I; |x-b\cdot c/d|\le\sqrt b/d^2\}
\ED
and the union
\BD
   \FF=\bigcup_{1\le d\le\sqrt b}\:\bigcup_{\genfrac{}{}{0pt}{1}{0\le c<d,}{(c,d)=1}}I_{c/d}.
\ED
Then for all integers $a$ in $I\sminus\FF$ with $(a,b)=1$, we have $|S(a,b)|\le 3\sqrt b+5$ (see \cite{Gi}). In other words, the Farey fractions $c/d$ of order $\lfloor \sqrt b\rfloor$
determine the intervals $I_{c/d}$ that contain all integers $a$ such that $|S(a,b)|$ is substantially larger than $\sqrt b$.
\DED sums inside and outside $\FF$ are illustrated by the diagrams in \cite{Gi}. The main tools used in the proof of this result are (\ref{2.3}) and a basic fact about Farey fractions.
Note that the number of integers inside $\FF$ is $\ll \sqrt b\log b$.

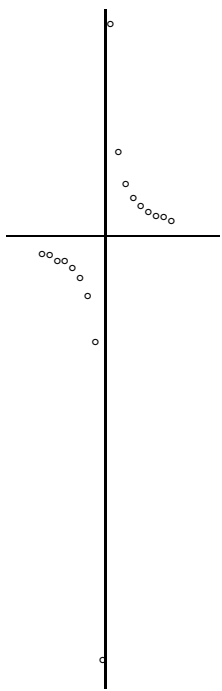
\unitlength1cm
\begin{picture}(12,10)
\put(2,6.5){\begin{picture}(8,8)
\put(-0.4,0){\line(1,0){2.8}}
\put(0.9,-6){\line(0,1){9}}
\put(3,0){$b=6761$:}
\put(3,-0.5){$(b,S(a,b))$ for $a$ inside $I_{2/3}$,}
\put(3,-1){$S(a,b)$ between $\approx -2251.33$ and $\approx 1123.17$,}
\put(3,-1.5){length of $I_{c/d}$: $\approx 18.27$}
\PT{  0.07}{ -0.24}
\PT{  0.17}{ -0.25}
\PT{  0.27}{ -0.33}
\PT{  0.37}{ -0.34}
\PT{  0.47}{ -0.43}
\PT{  0.57}{ -0.56}
\PT{  0.67}{ -0.80}
\PT{  0.77}{ -1.41}
\PT{  0.87}{ -5.63}
\PT{  0.97}{  2.81}
\PT{  1.07}{  1.12}
\PT{  1.17}{  0.69}
\PT{  1.27}{  0.50}
\PT{  1.37}{  0.40}
\PT{  1.47}{  0.32}
\PT{  1.57}{  0.26}
\PT{  1.67}{  0.25}
\PT{  1.77}{  0.19}
\end{picture}}

\put(0,0){\begin{picture}(0,0)
\end{picture}}
\end{picture}

The points $(a,S(a,b))$ with $a$ inside the interval $I_{c/d}$
can be described by means of the hyperbola
\BD
   H_{c/d}=\{(x,y); (x-b\cdot c/d)\cdot y=b/d^2\},
\ED
whose midpoint is $(b\cdot c/d,0)$.
Indeed, for $a\in I_{c/d}$, $a<b\cdot c/d$, $(a,S(a,b))$ lies close to the \BQ negative\EQ branch $\{(x,y)\in H_{c/d}; y<0\}$ and close to the vertical asymptote
$\{(b\cdot c/d,y); y\in \R\}$ of $H_{c/d}$. Conversely, $(a,S(a,b))$ lies close to the \BQ positive\EQ branch $\{(x,y)\in H_{c/d}; y>0\}$ and close to the same asymptote
for $a\in I_{c,d}$, $a>b\cdot c/d$. In general, the points $(a,S(a,b))$ are less close to the hyperbola if $a$ lies near the endpoints of $I_{c/d}$. These somewhat vague
statements have a more precise asymptotic meaning (see \cite{Gi}).

For small values of $d$ the course of the relevant parts of the hyperbola can be seen from the points $(a, S(a,b))$ with $a\in I_{c/d}$ (see the diagram, whose horizontal scale differs from the
vertical one).
But if $d>\sqrt 2\cdot b^{1/4}$, the interval $I_{c/d}$ contains at most one integer. Hence for most of the intervals $I_{c/d}$ there is at most one
point $(a,S(a,b))$ with $a\in I_{c/d}$ close to the hyperbola $H_{c/d}$.

In order to count large absolute values of \DED sums, we consider
\BD
  M_{\alpha}=\{a; 0\le a<b; (a,b)=1, |S(a,b)|>b^{\alpha}\}
\ED
for some $\alpha>0$. One can show
\BD
 |M_{\alpha}|\ge C_{\alpha}\phi(b)\log b/b^{\alpha}
\ED
for $1/3<\alpha<1$, if $b\to\infty$ (see \cite{Gi7}). This order of magnitude is in some sense best possible. It would be desirable to have an analogous result for $\alpha\le 1/3$,
but we only have
\BD
 |M_{\alpha}|\ge C_{\alpha}\phi(b)/b^{\alpha}
\ED
for $0<\alpha\le 1/3$. So the $\log$ factor has been lost in this result.

One of the few counting results about {\em small} values of \DED sums that we know is given in \cite{Va}. Its proof uses real-analytic modular forms.
For $X>0$, let
\BD
   A(X)=\{(a,b); 0\le a<b\le X, (a,b)=1\}.
\ED
Then for each $\alpha>0$
\BD
  \lim_{X\to\infty}\frac{|\{(a,b)\in A(X); S(a,b)<\alpha\log b\}|}{|A(X)|}=\frac 1{\pi}\cdot \arctan\left(\frac{\pi\alpha}6\right)+\frac 12.
\ED
This result implies
\BD
  \lim_{X\to\infty}\frac{|\{(a,b)\in A(X); |S(a,b)|<\alpha\log b\}|}{|A(X)|}=\frac 2{\pi}\cdot \arctan\left(\frac{\pi\alpha}6\right).
\ED
Computations suggest that a similar limiting behaviour takes place for {\em one} parameter $b$ instead of {\em all} $b<X$,
for instance,
\BD
 \lim_{b\to \infty}\frac{|\{a; 0\le a<b, (a,b)=1, |S(a,b)|<\log b\}|}{\phi(b)}=\frac 2{\pi}\cdot \arctan\left(\frac{\pi}6\right)=0.307072\LD
\ED
But this is far from being proved.

From (\ref{4.3}) we obtain, in a straightforward manner, the following result. Let $\alpha>0$. For every $\EPS>0$,
\BD
   \frac {|\{ (a,b)\in A(X); |S(a,b)|\le \alpha\}|}{|A(X)|}= \frac {\alpha}{3\log X}+E(\EPS,X)+o\left(\frac 1{\log X}\right)
\ED
with
\BD
  |E(\EPS,X)|\le \EPS/\log X,
\ED
as $X\to\infty$. It would be nice if the $E$-term could be omitted.

\section*{7. \DED sums near quadratic irrationals}

Let $\alpha\in\R$ be a quadratic irrational. For each $x\in\R$ there are values $S(a,b)$ arbitrarily close to $x$
such that $a/b$ is arbitrarily close to $\alpha$. This follows from the density of the set $\{(a/b,S(a,b)); a, b\in\Z, b>0, (a,b)=1\}$
in $\R^2$ (see Section 4). More interesting is the case when $a/b$ runs through the sequence of {\em best approximations} of $\alpha$,
i.e., the sequence of {\em convergents} of the continued fraction expansion of $\alpha$ (see \cite[p. 20 ff.]{RoSz}).  In this case it may happen that
the \DED sums $S(a,b)$ are concentrated near finitely many cluster points (see \cite {Gi8}). We consider only a special case here.

Let $\alpha$ be a quadratic irrational with continued fraction expansion $\alpha=[0,\OV{c_1,c_1,\LD,c_l}]$. So $c_1,\LD,c_l$ is the repeating block of $\alpha$. Let the convergents
$a_k/b_k$, $k=0,1,2,\LD$ be defined as usual, namely,
\BD
  a_{-1}=1,\,a_{0}=c_{0},\, b_{-1}=0,\, b_{0}=1,
\ED
and
\BD
  a_k=c_ka_{k-1}+a_{k-2},\, b_k=c_kb_{k-1}+b_{k-2}
\ED
for $k\ge 1$ (observe that $c_{l+1}=c_1, c_{l+2}=c_2$, and so on). The numbers $a_k,b_k$ are integers, $b_k\ge 1$, $(a_k,b_k)=1$.

Now suppose that $l$ is odd. We put $L=2l$. Then for each $j\in\{1,2,\LD,L\}$, the sequence $S(a_k,b_k)$
converges to
\BD
  \sum_{r=1}^j(-1)^{r-1}c_r+\alpha+ \begin{cases} 1/\alpha_{j}-3, & \MB{ if } j \MB{ is odd};\\
                                                 -1/\alpha_{j}, & \MB { otherwise,}
                                    \end{cases}
\ED
for $k\to\infty$, $k\equiv j\mod L$. Here $\alpha_j$ is the quadratic irrational
\BD
 \alpha_j=[c_j,c_{j-1},\LD,c_1,\OV{c_l,c_{l-1},\LD,c_{1}}].
\ED
This follows from formula (\ref{2.5}).
Accordingly, the sequence $S(a_k,b_k)$ is bounded and has at most $L$ different cluster points when $k$ tends to infinity.

\MN
{\em Example.} Let $\alpha=[0,\OV{1,2,2}]=(\sqrt{85}-5)/6=0.703257\LD$ In this case the six cluster points are
$(\sqrt{85}-11)/3=-0.5934851\LD$, -2/3, $(4\sqrt{85}-50)/15=-0.8747881\LD$, $0$, $\sqrt{85}/3-3=0.0731848\LD$,
$(\sqrt{85}-5)/15=0.2813029\LD$ The \DED sums $S(a_k,b_k)$, $k=7,\LD,12$, are already quite close to these cluster points.


\vspace{0.5cm}
\noindent
Kurt Girstmair            \\
Institut f\"ur Mathematik \\
Universit\"at Innsbruck   \\
Technikerstr. 13/7        \\
A-6020 Innsbruck, Austria \\
Kurt.Girstmair@uibk.ac.at

\end{document}